\documentclass{amsart}
\usepackage{amssymb}

\usepackage{amsmath, amssymb, euscript,  enumerate}

\newtheorem{theorem}{Theorem}
\newtheorem{lemma}{Lemma}
\newtheorem{corollary}{Corollary}

\begin{document}
\author{Ovidiu Munteanu and Mu-Tao Wang}
\title{The curvature of gradient Ricci solitons}
\date{}

\begin{abstract}
We study integral and pointwise bounds on the curvature of gradient
shrinking Ricci solitons. As applications, we discuss gap and compactness
results for gradient shrinkers.
\end{abstract}

\maketitle

\section{Introduction}

Our goal in this paper is to obtain further information about the curvature
of gradient shrinking Ricci solitons. This is important for a better
understanding and ultimately for the classification of these manifolds. The
classification of gradient shrinkers is known in dimensions 2 and 3, and
assuming locally conformally flatness, in all dimensions $n\geq 4$ (see \cite
{P, NW, CWZ, PW, Zh, MS, M}). Many of the techniques used in these works
required some control of the Ricci curvature. For example, in \cite{PW}
gradient shrinking Ricci solitons which are locally conformally flat were
classified assuming an integral condition on the Ricci tensor. This
condition and other integral estimates of the curvature were later proved in 
\cite{MS}.

Without making the strong assumption of being conformally flat, it is
natural to ask whether similar estimates are true for the Riemann curvature
tensor. In this paper we are able to prove pointwise estimates on the
Riemann curvature, assuming in addition that the Ricci curvature is bounded.
We will show that any gradient shrinking Ricci soliton with bounded Ricci
curvature has Riemann curvature tensor growing at most polynomially in the
distance function. We note that by Shi's local derivative estimates we can
then obtain growth estimates on all derivatives of the curvature. This, in
particular, proves weighted $L^2$ estimates for the Riemann curvature tensor
and its covariant derivatives.

We point out that for self shrinkers of the mean curvature flow Colding and
Minicozzi \cite{CM} were able to prove weighted $L^2$ estimates for the
second fundamental form, assuming the mean curvature is positive. These
estimates were instrumental in the classification of stable shrinkers. Our
estimates can be viewed as parallel to theirs, however the classification of
gradient Ricci solitons is still a major open question in the field.

A gradient shrinking Ricci soliton is a Riemannian manifold $\left(
M,g\right) $ for which there exists a potential function $f$ such that 
\begin{equation}  \label{soliton_eq}
R_{ij}+f_{ij}=\frac 12g_{ij}.
\end{equation}

It can be shown directly from the equation that the quantity $\left| \nabla
f\right| ^2+R-f$ is constant on $M$ hence we can normalize $f$ such that 
\begin{equation*}
\left| \nabla f\right| ^2+R=f.
\end{equation*}

Let us denote with $\left| Rc\right| $ and $\left| Rm\right| $ the norms of
the Ricci and Riemann tensors, respectively. We now state the main result of
this paper.

\begin{theorem}
Let $\left( M,g,f\right) $ be a gradient shrinking Ricci soliton with
bounded Ricci curvature. Then the Riemann curvature tensor grows at most
polynomially in the distance function i.e. 
\begin{equation*}
\left| Rm\right| \left( x\right) \leq C\left( r\left( x\right) +1\right) ^a,
\end{equation*}
for some constant $a>0.$
\end{theorem}

We apply this theorem to prove a gap result for gradient shrinkers. We show
that if the Ricci curvature is small enough everywhere on $M$ then the
soliton is isometric to the Gaussian soliton $\left( \mathbb{R}^n,dx^2,\frac
14\left| x\right| ^2\right) .$

\begin{corollary}
Let $\left( M,g,f\right) $ be a gradient shrinking Ricci soliton. If $\left|
Rc\right| \leq \frac 1{100n^2}$ on $M$ then $M$ is isometric to the Gaussian
soliton.
\end{corollary}

Yokota \cite{Y} has obtained a gap theorem for gradient shrinking Ricci
solitons, namely he showed that if $\left( 4\pi \right) ^{-\frac
n2}\int_Me^{-f}>1-\varepsilon _n$ for some $\varepsilon _n$ depending on $n$
then $M$ is isometric to the Gaussian soliton.

Another application of our main theorem is in relation to compactness
results for Ricci solitons. This topic has been studied recently in both
compact (\cite{CS, W, Z}) and noncompact (\cite{HM}) settings. We recall a
recent result for complete noncompact shrinkers, due to Haslhofer and
M\"{u}ller \cite{HM}. Let $\left( M_i,g_i,\bar{f}_i\right) $ a sequence of
gradient shrinkers with the potentials $\bar{f}_i$ normalized such that $%
\left( 4\pi \right) ^{-\frac n2}\int_{M_i}e^{-\bar{f}_i}=1$. Consider $z_i$
a point where $\bar{f}_i$ attains its minimum. Assume that Perelman's
entropy $\mu _i$ has a uniform bound from below i.e. there exists $\bar{\mu}$
such that 
\begin{equation*}
\mu _i:=\left( 4\pi \right) ^{-\frac n2}\int_{M_i}\left( \left| \nabla \bar{f%
}_i\right| ^2+R_{g_i}+\bar{f}_i-n\right) e^{-\bar{f}_i}\geq \bar{\mu}.
\end{equation*}
If, moreover, for any $i$ and $r>0$ 
\begin{equation*}
\int_{B_{z_i}\left( r\right) }\left| Rm_{g_i}\right| ^{n/2}\leq E\left(
r\right) ,
\end{equation*}
then a subsequence of $\left( M_i,g_i,\bar{f}_i,z_i\right) $ converges to an
orbifold gradient shrinker in the pointed Cheeger-Gromov sense.

As a consequence of our main theorem, we can show the following result. Here
it is not important how we normalize the function $f,$ but we need to take $%
x_0$ a point where $f$ assumes its minimum value on $M.$ Moreover, here we
also need to take $n\geq 6.$ Notice however that in dimensions 2 and 3
shrinking solitons are completely classified and in dimension 4, Haslhofer
and M\"{u}ller have proved that $\int_{B_{z_i}\left( r\right) }\left|
Rm_{g_i}\right| ^2\leq E\left( r\right) $ using the Gauss-Bonnet theorem.
However, it seems that their argument is more special and works only in
lower dimensions.

\begin{corollary}
Let $\left( M,g,f\right) $ be a gradient shrinker with $dimM=n\geq 6$ and $%
\left| Rc\right| \leq K$ on $M$. There exists a $r_0$ depending on $n$ and $K
$ such that if $\int_{B_{x_0}\left( r_0\right) }\left| Rm\right| ^{n/2}\leq L
$ for a minimum point $x_0$ of $f$, then for any $r>0$ we have 
\begin{equation*}
\int_{B_{x_0}\left( r\right) }\left| Rm\right| ^{n/2}\leq E\left( r\right) ,
\end{equation*}
where $E$ depends on $n,K$ and $L$.
\end{corollary}

In particular, this Corollary and the main theorem in \cite{HM} implies
compactness of shrinkers assuming Ricci curvature bounds and only local
curvature bounds.

\begin{corollary}
Let $\left( M_i,g_i,\bar{f}_i\right) $ be a sequence of gradient shrinking
Ricci solitons normalized by $\left( 4\pi \right) ^{-\frac n2}\int_{M_i}e^{-%
\bar{f}_i}=1.$ Assume $\mu _i\geq \bar{\mu}$ and that we have a uniform
bound on the Ricci curvature $\left| Rc_{g_i}\right| \leq K$. Then there
exists $r_0$ depending on $n$ and $K$ such that if 
\begin{equation*}
\sup_i\int_{B_{z_i}\left( r_0\right) }\left| Rm_{g_i}\right| ^{n/2}
\end{equation*}
is finite for $z_i$ a minimum point of $\bar{f}_i$, then a subsequence of $%
\left( M_i,g_i,\bar{f}_i,z_i\right) $ converges to an orbifold gradient
shrinker in the pointed Cheeger-Gromov sense.
\end{corollary}

\section{Proof of the curvature estimate}

The idea of the proof of Theorem 1 is the following.

From the Ricci soliton equation \eqref{soliton_eq}, we can estimate $\Delta
_f\left| Rm\right| ^2\geq -c\left| Rm\right| ^3,$ where $\Delta _f=\Delta
-\nabla f\cdot \nabla .$ It is natural therefore to attempt to use Moser
iteration for this problem. Since Ricci is bounded below, the Sobolev
constant is uniformly bounded on arbitrary balls of fixed radius $=1.$
Therefore it is known that the Moser iteration will work if we can control
the $L^p$ norm (for $p>n/2$) of $\left| Rm\right| $ on any ball of radius
one. This is quite technical and it is done in the second Lemma below.

To get $L^p$ estimates, we use the Ricci soliton equation and that the Ricci
curvature is bounded. At the core of our estimates is a formula that relates
the divergence of the Riemann curvature tensor of a Ricci soliton to the
gradient of its Ricci curvature.

We first prove some Lemmas. 
Everywhere in this section $\left( M,g,f\right) $
is a gradient shrinking Ricci soliton with bounded Ricci curvature. 
We use
the notation: 
\begin{eqnarray*}
\left| Rc\right| ^2 &=&\sum \left| R_{ij}\right| ^2, \\
\left| \nabla Rc\right| ^2 &=&\sum \left| \nabla _kR_{ij}\right| ^2, \\
\left| Rm\right| ^2 &=&\sum \left| R_{ijkl}\right| ^2, \\
\left| \nabla Rm\right| ^2 &=&\sum \left| \nabla _hR_{ijkl}\right| ^2.
\end{eqnarray*}

We recall some basic identities for shrinking Ricci solitons, which are
essential in the proof. First, taking trace of the soliton equation we get $%
R+\Delta f=\frac n2.$ As we have mentioned above, using the Bianchi
identities and normalizing $f$ we get $\left| \nabla f\right| ^2+R=f.$ This
normalization of $f$ will be assumed throughout the paper. Other formulas
that follow from Bianchi and Ricci identities and the soliton equation are 
\begin{eqnarray*}
\nabla _iR &=&2R_{ij}f_j, \\
\nabla _kR_{jk} &=&R_{jk}f_k, \\
\nabla _lR_{ijkl} &=&R_{ijkl}f_l, \\
\nabla _jR_{ki}-\nabla _iR_{kj} &=&R_{ijkl}f_l.
\end{eqnarray*}
We will not prove these here since they are quite standard, see e.g. \cite
{ENM}.

\begin{lemma}
We have: 
\begin{eqnarray}
\left| \nabla Rc\right| ^2 &\leq &\frac 12\Delta \left| Rc\right| ^2-\frac
12\nabla f\cdot \nabla \left| Rc\right| ^2+c\left| Rm\right|   \label{f1} \\
\left| \nabla Rm\right| ^2 &\leq &\frac 12\Delta \left| Rm\right| ^2-\frac
12\nabla f\cdot \nabla \left| Rm\right| ^2+c\left| Rm\right| ^3  \label{f2}
\end{eqnarray}
\end{lemma}

Proof. It is known (see \cite{PW}) that 
\begin{eqnarray*}
\Delta R_{ij} &=&\nabla f\cdot \nabla R_{ij}+R_{ij}-2R_{ikjh}R_{kh} \\
\Delta Rm &=&\nabla f\cdot \nabla Rm+Rm-2\left( Rm^2+Rm^{\#}\right) .
\end{eqnarray*}
The Lemma follows immediately from here. $\blacksquare $

We now prove the following estimate, which is of independent interest. 

\begin{lemma}
For any $p\geq 2$ there exist positive constants $C$ and $a$ such that 
\begin{equation*}
\int_M\left| Rm\right| ^p\left( f+1\right) ^{-a}\leq C.
\end{equation*}
In particular, there exist positive constants $C$ and $a$ such that for any $%
x\in M$ we have: 
\begin{equation*}
\int_{B_x\left( 1\right) }\left| Rm\right| ^p\leq C\left( r\left( x\right)
+1\right) ^{2a}.
\end{equation*}
\end{lemma}

Proof. Let us denote by $\rho :=2\sqrt{f}$ and for $r>>1$ let 
\begin{equation*}
D\left( r\right) :=\left\{ x\in M:\;\rho \left( x\right) \leq r\right\} .
\end{equation*}
Notice that $D\left( r\right) $ is always compact, in fact approximates well
the geodesic ball of radius $r$ when $r$ is large. Here we recall that $f$
has the following asymptotics, see \cite{CZ}: 
\begin{equation*}
\left( \frac 12r\left( x\right) -c\right) ^2\leq f\left( x\right) \leq
\left( \frac 12r\left( x\right) +c\right) ^2,\;\;\text{for}\;\;r\left(
x\right) \geq r_0.
\end{equation*}
Here and below we denote by $r\left( x\right) $ the distance from $x$ to a
fixed point $x_0\in M.$ Let us point out moreover that if $x_0$ is chosen to
be a minimum point of $f$ then $c$ and $r_0$ will depend only on $n,$ see 
\cite{HM}.

We take the following cut-off 
\begin{equation*}
\phi =\left\{ 
\begin{array}{c}
\frac 1{r^2}\left( \frac 14r^2-f\left( x\right) \right) \\ 
0
\end{array}
\right. 
\begin{array}{l}
\text{if} \\ 
\text{if}
\end{array}
\begin{array}{l}
x\in D\left( r\right) \\ 
x\in M\backslash D\left( r\right)
\end{array}
\end{equation*}

Since the Ricci curvature is bounded on $M$, let us set 
\begin{equation*}
K=\sup_M\left| Rc\right| .
\end{equation*}

Let $a$ be a fixed number to be determined later, depending on $n,K$ and $p.$
Consider also $q$ a large enough integer, $q\geq 2p+1$. We will discuss
first the case when $p\geq 3,$ the case $2\leq p<3$ will follow immediately
by H\"{o}lder's inequality.

We have, integrating by parts, that 
\begin{gather}
a\int_M\left| Rm\right| ^p\left| \nabla f\right| ^2\left( f+1\right)
^{-a-1}\phi ^q=-\int_M\left| Rm\right| ^p\nabla f\cdot \nabla \left(
f+1\right) ^{-a}\phi ^q  \label{a} \\
=\int_M\left| Rm\right| ^p\left( \Delta f\right) \left( f+1\right) ^{-a}\phi
^q+\int_M\left| Rm\right| ^p\left( f+1\right) ^{-a}\nabla f\cdot \nabla \phi
^q  \notag \\
+\int_M\left( \nabla \left| Rm\right| ^p\cdot \nabla f\right) \left(
f+1\right) ^{-a}\phi ^q.  \notag
\end{gather}
Let us explain why we take these functions in (\ref{a}). We take $\left|
\nabla f\right| ^2$ because we want to use integration by parts and the
symmetries of the Riemann curvature tensor (which are implied by the soliton
equation). For shrinking solitons the factor $e^{-f}$ seems much more
convenient than $\left( f+1\right) ^{-a},$ however the former factor only
gives exponential growth control. So to prove weighted $L^2$ estimates for
the Riemann curvature it is easier to start from $\int_M\left| Rm\right|
^p\left| \nabla f\right| ^2e^{-f}\phi ^q$ and carry all the estimates below,
with some useful simplifications, but in order to prove Lemma 2 we need to
use the less natural weight $\left( f+1\right) ^{-a},$ at the expense of
more complicated computations.

Everywhere in this section $c$ will denote a constant that depends on $n$, $%
p $, $q$ and $K$ but not on $a$. We will use $C_1,$ $C_2,...$ etc. to denote
finite constants that have a more complicated dependence, such as on $%
\sup_\Omega \left| Rm\right| $ over some compact set $\Omega .$ However, we
stress that all the constants $c$ or $C_1,C_2...$ are independent of $r.$

We now check that: 
\begin{equation*}
a\left| \nabla f\right| ^2\left( f+1\right) ^{-a-1}-\Delta f\left(
f+1\right) ^{-a}=\left( a\frac{f-R}{f+1}-\left( \frac n2-R\right) \right)
\left( f+1\right) ^{-a}.
\end{equation*}
Furthermore, since $R\geq 0$ on any gradient Ricci soliton ( \cite{C, Ca})
we see that there exists a constant $r_1$ depending on $n,K$ and $a$ (e.g., $%
r_1=\sqrt{8a\left( K+1\right) }$ ) such that on $M\backslash D\left(
r_1\right) $ we have 
\begin{equation*}
a\left| \nabla f\right| ^2\left( f+1\right) ^{-a-1}-\Delta f\left(
f+1\right) ^{-a}\geq \left( a-n\right) \left( f+1\right) ^{-a}.
\end{equation*}

Notice also that by the choice of cut-off we have: 
\begin{equation*}
\nabla f\cdot \nabla \phi ^q=-\frac q{r^2}\phi ^{q-1}\left| \nabla f\right|
^2\leq 0.
\end{equation*}
Using these simple estimates in (\ref{a}) it follows that 
\begin{equation}
\left( a-n\right) \int_M\left| Rm\right| ^p\left( f+1\right) ^{-a}\phi
^q\leq \int_M\left( \nabla \left| Rm\right| ^p\cdot \nabla f\right) \;\left(
f+1\right) ^{-a}\phi ^q+C_1.  \label{b}
\end{equation}
Here we have set 
\begin{equation*}
C_1:=\int_{D\left( r_1\right) }\left( -a\left| \nabla f\right| ^2\left(
f+1\right) ^{-1}+\Delta f+a-n\right) \left| Rm\right| ^p\left( f+1\right)
^{-a}\phi ^q.
\end{equation*}

Let us compute, using the Bianchi identities: 
\begin{eqnarray*}
\nabla \left| Rm\right| ^2\cdot \nabla f &=&2f_h\left( \nabla
_hR_{ijkl}\right) R_{ijkl} \\
&=&4f_h\left( \nabla _lR_{ijkh}\right) R_{ijkl}
\end{eqnarray*}
Therefore, using this in the right hand side of (\ref{b}) we get 
\begin{gather*}
\int_M\left( \nabla \left| Rm\right| ^p\cdot \nabla f\right) \;\left(
f+1\right) ^{-a}\phi ^q \\
=2p\int_Mf_h\left( \nabla _lR_{ijkh}\right) R_{ijkl}\left| Rm\right|
^{p-2}\left( f+1\right) ^{-a}\phi ^q \\
=-2p\int_MR_{ijkh}\nabla _l\left( f_hR_{ijkl}\left| Rm\right| ^{p-2}\left(
f+1\right) ^{-a}\phi ^q\right)  \\
=-2p\int_MR_{ijkh}f_{hl}R_{ijkl}\left| Rm\right| ^{p-2}\left( f+1\right)
^{-a}\phi ^q \\
-2p\int_MR_{ijkh}f_h\left( \nabla _lR_{ijkl}\right) \left| Rm\right|
^{p-2}\left( f+1\right) ^{-a}\phi ^q \\
-2p\int_MR_{ijkh}f_hR_{ijkl}\left( \nabla _l\left| Rm\right| ^{p-2}\right)
\left( f+1\right) ^{-a}\phi ^q \\
+2ap\int_MR_{ijkh}f_hR_{ijkl}f_l\left| Rm\right| ^{p-2}\left( f+1\right)
^{-a-1}\phi ^q \\
-2pq\int_MR_{ijkh}f_hR_{ijkl}\phi _l\left| Rm\right| ^{p-2}\left( f+1\right)
^{-a}\phi ^{q-1} \\
=I+II+III+IV+V.
\end{gather*}
It is easy to see that since the Ricci curvature is bounded, 
\begin{equation*}
I=-2p\int_MR_{ijkh}f_{hl}R_{ijkl}\left| Rm\right| ^{p-2}\left( f+1\right)
^{-a}\phi ^q\leq c\int_M\left| Rm\right| ^p\left( f+1\right) ^{-a}\phi ^q.
\end{equation*}
Furthermore, using that for a gradient shrinker we have (see \cite{ENM, MS}) 
\begin{equation*}
\nabla _lR_{ijkl}=R_{ijkl}f_l
\end{equation*}
we see that 
\begin{equation*}
II+IV=-2p\int_MR_{ijkh}f_hR_{ijkl}f_l\left| Rm\right| ^{p-2}\left( 1-\frac
a{f+1}\right) \left( f+1\right) ^{-a}\phi ^q\leq C_2,
\end{equation*}
where we have set 
\begin{equation*}
C_2:=2p\int_{D\left( 2\sqrt{a-1}\right) }\left| R_{ijkh}f_h\right| ^2\left|
Rm\right| ^{p-2}\left( \frac a{f+1}-1\right) \left( f+1\right) ^{-a}\phi ^q.
\end{equation*}
The estimate above follows because if $f\left( x\right) +1>a$ then the
integral is negative. Let us point out however that in fact each of $\left|
II\right| $ and $\left| IV\right| $ can be estimated by a similar argument
as in the proof of inequality (\ref{f}) below. Clearly, since $a$ is
independent of $r$, so are the constants $C_1$ and $C_2$ obtained so far.

We use the above estimates in (\ref{b}), and get that 
\begin{gather}
\left( a-c\right) \int_M\left| Rm\right| ^p\left( f+1\right) ^{-a}\phi ^q
\label{c} \\
\leq -2p\int_MR_{ijkh}f_hR_{ijkl}\left( \nabla _l\left| Rm\right|
^{p-2}\right) \left( f+1\right) ^{-a}\phi ^q  \notag \\
+\frac{2pq}{r^2}\int_M\left| R_{ijkh}f_h\right| ^2\left| Rm\right|
^{p-2}\left( f+1\right) ^{-a}\phi ^{q-1}+C_1+C_2.  \notag
\end{gather}
Recall that for Ricci solitons we have (see e.g. \cite{ENM, MS}) 
\begin{equation*}
R_{ijkh}f_h=\nabla _jR_{ik}-\nabla _iR_{kj},
\end{equation*}
so that we can estimate the first term in (\ref{c}) by 
\begin{gather}
-2p\int_MR_{ijkh}f_hR_{ijkl}\left( \nabla _l\left| Rm\right| ^{p-2}\right)
\left( f+1\right) ^{-a}\phi ^q  \label{d} \\
\leq c\int_M\left| \nabla Rc\right| \left| \nabla Rm\right| \left| Rm\right|
^{p-2}\left( f+1\right) ^{-a}\phi ^q  \notag \\
\leq c\int_M\left| \nabla Rc\right| ^2\left| Rm\right| ^{p-1}\left(
f+1\right) ^{-a}\phi ^q+c\int_M\left| \nabla Rm\right| ^2\left| Rm\right|
^{p-3}\left( f+1\right) ^{-a}\phi ^q.  \notag
\end{gather}
We now work on the second term in (\ref{c}), and we will use below some
interpolations which are important throughout the rest of the proof.

\begin{gather}
\frac 1{r^2}\int_MR_{ijkh}f_hR_{ijkl}f_l\left| Rm\right| ^{p-2}\left(
f+1\right) ^{-a}\phi ^{q-1}  \notag \\
=\frac 2{r^2}\int_M\left( \nabla _jR_{ik}\right) R_{ijkl}f_l\left| Rm\right|
^{p-2}\left( f+1\right) ^{-a}\phi ^{q-1}  \notag \\
=-\frac 2{r^2}\int_MR_{ik}\nabla _j\left( R_{ijkl}f_l\left| Rm\right|
^{p-2}\left( f+1\right) ^{-a}\phi ^{q-1}\right)  \notag \\
=-\frac 2{r^2}\int_MR_{ik}f_{lj}R_{ijkl}\left| Rm\right| ^{p-2}\left(
f+1\right) ^{-a}\phi ^{q-1}  \notag \\
-\frac 2{r^2}\int_MR_{ik}f_l\left( \nabla _jR_{ijkl}\right) \left| Rm\right|
^{p-2}\left( f+1\right) ^{-a}\phi ^{q-1}  \notag \\
-\frac 2{r^2}\int_MR_{ik}R_{ijkl}f_l\left( \nabla _j\left| Rm\right|
^{p-2}\right) \left( f+1\right) ^{-a}\phi ^{q-1}  \notag \\
+\frac{2a}{r^2}\int_MR_{ik}R_{ijkl}f_lf_j\left| Rm\right| ^{p-2}\left(
f+1\right) ^{-a-1}\phi ^{q-1}  \notag \\
-\frac{2(q-1)}{r^2}\int_MR_{ik}R_{ijkl}f_l\phi _j\left| Rm\right|
^{p-2}\left( f+1\right) ^{-a}\phi ^{q-2}  \notag \\
=I+II+III+IV+V.  \label{e}
\end{gather}
Since the Ricci curvature is bounded, we get 
\begin{eqnarray*}
I &=&-\frac 2{r^2}\int_MR_{ik}f_{lj}R_{ijkl}\left| Rm\right| ^{p-2}\left(
f+1\right) ^{-a}\phi ^{q-1} \\
&\leq &\frac c{r^2}\int_M\left| Rm\right| ^{p-1}\left( f+1\right) ^{-a}\phi
^{q-1}.
\end{eqnarray*}
Using again that $\nabla _lR_{ijkl}=R_{ijkl}f_l$ and that $\left| \nabla
f\right| ^2\leq f\leq \frac 14r^2$ on $D\left( r\right) ,$ it follows 
\begin{equation*}
II\leq c\int_M\left| Rm\right| ^{p-1}\left( f+1\right) ^{-a}\phi ^{q-1}.
\end{equation*}
Furthermore, 
\begin{eqnarray*}
III &=&-\frac 2{r^2}\int_MR_{ik}R_{ijkl}f_l\left( \nabla _j\left| Rm\right|
^{p-2}\right) \left( f+1\right) ^{-a}\phi ^{q-1} \\
&\leq &c\int_M\left| \nabla Rm\right| \left| Rm\right| ^{p-2}\left(
f+1\right) ^{-a}\phi ^{q-1} \\
&\leq &c\int_M\left| \nabla Rm\right| ^2\left| Rm\right| ^{p-3}\left(
f+1\right) ^{-a}\phi ^q+c\int_M\left| Rm\right| ^{p-1}\left( f+1\right)
^{-a}\phi ^{q-2}.
\end{eqnarray*}
Similarly, 
\begin{eqnarray*}
IV &=&\frac{2a}{r^2}\int_MR_{ik}R_{ijkl}f_lf_j\left| Rm\right| ^{p-2}\left(
f+1\right) ^{-a-1}\phi ^{q-1} \\
&\leq &\frac{ca}{r^2}\int_M\left| Rm\right| ^{p-1}\left( f+1\right)
^{-a}\phi ^{q-1} \\
&\leq &c\int_M\left| Rm\right| ^{p-1}\left( f+1\right) ^{-a}\phi ^{q-1},
\end{eqnarray*}
by taking $r\geq \sqrt{a}.$ Finally, we also have 
\begin{equation*}
V\leq \frac c{r^2}\int_M\left| Rm\right| ^{p-1}\left( f+1\right) ^{-a}\phi
^{q-2}.
\end{equation*}
Using these estimates in (\ref{e}) we get 
\begin{gather}
\frac{2pq}{r^2}\int_M\left| R_{ijkh}f_h\right| ^2\left| Rm\right|
^{p-2}\left( f+1\right) ^{-a}\phi ^{q-1}  \label{f} \\
\leq c\int_M\left| \nabla Rm\right| ^2\left| Rm\right| ^{p-3}\left(
f+1\right) ^{-a}\phi ^q  \notag \\
+c\int_M\left| Rm\right| ^{p-1}\left( f+1\right) ^{-a}\phi ^{q-2}.  \notag
\end{gather}
Notice moreover that we can interpolate, using Young's inequality: 
\begin{gather}
\int_M\left| Rm\right| ^{p-1}\left( f+1\right) ^{-a}\phi ^{q-2}\leq
\int_M\left| Rm\right| ^{p-1}\phi ^{q\frac{p-1}p}\phi ^{\frac qp-2}\left(
f+1\right) ^{-a}  \notag \\
\leq \varepsilon \int_M\left| Rm\right| ^p\left( f+1\right) ^{-a}\phi
^q+c\left( \varepsilon \right) \int_M\left( f+1\right) ^{-a}\phi ^{q-2p}.
\label{g}
\end{gather}
We want to use this in (\ref{f}), hence here we can take $\varepsilon =1.$
Moreover, to guarantee that $\phi ^{q-2p}$ is well defined, we take $q\geq
2p+1$.

Plug (\ref{g}) in (\ref{f}) and then use (\ref{f}) and (\ref{d}) in (\ref{c}%
); it results that

\begin{gather}
\left( a-c\right) \int_M\left| Rm\right| ^p\left( f+1\right) ^{-a}\phi
^q\leq c\int_M\left| \nabla Rc\right| ^2\left| Rm\right| ^{p-1}\left(
f+1\right) ^{-a}\phi ^q  \label{h} \\
+c\int_M\left| \nabla Rm\right| ^2\left| Rm\right| ^{p-3}\left( f+1\right)
^{-a}\phi ^q+C_1+C_2+C_0.  \notag
\end{gather}
We have denoted with 
\begin{equation*}
C_0:=c\int_M\left( f+1\right) ^{-a},
\end{equation*}
and observe that taking $q=2p+1$ and taking $a$ such that $a>\frac n2+1$
then $C_0$ is a finite constant, independent of $r$. Indeed it is known that
the volume growth of $M$ is polynomial, see \cite{CZ}.

We finish the proof by estimating each of the two terms in the right hand
side of (\ref{h}). Start with the first, which by (\ref{f1}) we have: 
\begin{gather}
2\int_M\left| \nabla Rc\right| ^2\left| Rm\right| ^{p-1}\left( f+1\right)
^{-a}\phi ^q\leq \int_M\left( \Delta \left| Rc\right| ^2\right) \left|
Rm\right| ^{p-1}\left( f+1\right) ^{-a}\phi ^q  \label{i} \\
-\int_M\left( \nabla f\cdot \nabla \left| Rc\right| ^2\right) \left|
Rm\right| ^{p-1}\left( f+1\right) ^{-a}\phi ^q+c\int_M\left| Rm\right|
^p\left( f+1\right) ^{-a}\phi ^q  \notag
\end{gather}
Let us observe that 
\begin{gather*}
\int_M\left( \Delta \left| Rc\right| ^2\right) \left| Rm\right| ^{p-1}\left(
f+1\right) ^{-a}\phi ^q \\
=-\int_M\nabla \left| Rc\right| ^2\cdot \nabla \left( \left| Rm\right|
^{p-1}\left( f+1\right) ^{-a}\phi ^q\right) \\
=-\int_M\nabla \left| Rc\right| ^2\cdot \nabla \left| Rm\right| ^{p-1}\left(
f+1\right) ^{-a}\phi ^q \\
+a\int_M\left( \nabla \left| Rc\right| ^2\cdot \nabla f\right) \left|
Rm\right| ^{p-1}\left( f+1\right) ^{-a-1}\phi ^q \\
-q\int_M\left( \nabla \left| Rc\right| ^2\cdot \nabla \phi \right) \left|
Rm\right| ^{p-1}\left( f+1\right) ^{-a}\phi ^{q-1} \\
\leq c\int_M\left| \nabla Rc\right| \left| \nabla Rm\right| \left| Rm\right|
^{p-2}\left( f+1\right) ^{-a}\phi ^q \\
+ca\int_M\left| \nabla Rc\right| \left| Rm\right| ^{p-1}\left( f+1\right)
^{-a}\phi ^q \\
+\frac cr\int_M\left| \nabla Rc\right| \left| Rm\right| ^{p-1}\left(
f+1\right) ^{-a}\phi ^{q-1} \\
\leq c\int_M\left| \nabla Rc\right| \left| \nabla Rm\right| \left| Rm\right|
^{p-2}\left( f+1\right) ^{-a}\phi ^q \\
+ca\int_M\left| \nabla Rc\right| \left| Rm\right| ^{p-1}\left( f+1\right)
^{-a}\phi ^{q-1}.
\end{gather*}
Furthermore, let us use that 
\begin{eqnarray*}
&&c\int_M\left| \nabla Rc\right| \left| \nabla Rm\right| \left| Rm\right|
^{p-2}\left( f+1\right) ^{-a}\phi ^q \\
&\leq &\frac 14\int_M\left| \nabla Rc\right| ^2\left| Rm\right| ^{p-1}\left(
f+1\right) ^{-a}\phi ^q \\
&&+c\int_M\left| \nabla Rm\right| ^2\left| Rm\right| ^{p-3}\left( f+1\right)
^{-a}\phi ^q,
\end{eqnarray*}
and, similarly, 
\begin{eqnarray*}
&&ca\int_M\left| \nabla Rc\right| \left| Rm\right| ^{p-1}\left( f+1\right)
^{-a}\phi ^{q-1} \\
&\leq &\frac 14\int_M\left| \nabla Rc\right| ^2\left| Rm\right| ^{p-1}\left(
f+1\right) ^{-a}\phi ^q \\
&&+ca^2\int_M\left| Rm\right| ^{p-1}\left( f+1\right) ^{-a}\phi ^{q-2}.
\end{eqnarray*}
We conclude from above that 
\begin{eqnarray*}
&&\int_M\left( \Delta \left| Rc\right| ^2\right) \left| Rm\right|
^{p-1}\left( f+1\right) ^{-a}\phi ^q \\
&\leq &\frac 12\int_M\left| \nabla Rc\right| ^2\left| Rm\right| ^{p-1}\left(
f+1\right) ^{-a}\phi ^q+c\int_M\left| \nabla Rm\right| ^2\left| Rm\right|
^{p-3}\left( f+1\right) ^{-a}\phi ^q \\
&&+ca^2\int_M\left| Rm\right| ^{p-1}\left( f+1\right) ^{-a}\phi ^{q-2}.
\end{eqnarray*}
Moreover,

\begin{gather*}
-\int_M\left( \nabla f\cdot \nabla \left| Rc\right| ^2\right) \left|
Rm\right| ^{p-1}\left( f+1\right) ^{-a}\phi ^q \\
\leq c\int_M\left| \nabla Rc\right| \left| Rm\right| ^{p-1}\left| \nabla
f\right| \left( f+1\right) ^{-a}\phi ^q \\
\leq \frac 12\int_M\left| \nabla Rc\right| ^2\left| Rm\right| ^{p-1}\left(
f+1\right) ^{-a}\phi ^q+c\int_M\left| Rm\right| ^{p-1}\left| \nabla f\right|
^2\left( f+1\right) ^{-a}\phi ^q.
\end{gather*}
We can use the same idea as in (\ref{g}) to bound: 
\begin{eqnarray*}
&&\int_M\left| Rm\right| ^{p-1}\left| \nabla f\right| ^2\left( f+1\right)
^{-a}\phi ^q \\
&\leq &c\int_M\left| Rm\right| ^p\left( f+1\right) ^{-a}\phi
^q+c\int_M\left( f+1\right) ^{-a}\left| \nabla f\right| ^{2p}\phi ^q.
\end{eqnarray*}
Once again, we can take $a$ large enough e.g. $a>\frac n2+p+1$ so that 
\begin{equation*}
\int_M\left( f+1\right) ^{-a}\left| \nabla f\right| ^{2p}<\int_M\left(
f+1\right) ^{-a+p}<\infty .
\end{equation*}
We also use (\ref{g}) for $\varepsilon =\frac 1{a^2}$ to get 
\begin{eqnarray*}
&&ca^2\int_M\left| Rm\right| ^{p-1}\left( f+1\right) ^{-a}\phi ^{q-2} \\
&\leq &c\int_M\left| Rm\right| ^p\left( f+1\right) ^{-a}\phi
^q+ca^{2p}\int_M\left( f+1\right) ^{-a}\phi ^{q-2p}.
\end{eqnarray*}
Therefore, plugging all in (\ref{i}) gives 
\begin{gather*}
\int_M\left| \nabla Rc\right| ^2\left| Rm\right| ^{p-1}\left( f+1\right)
^{-a}\phi ^q\leq c\int_M\left| \nabla Rm\right| ^2\left| Rm\right|
^{p-3}\left( f+1\right) ^{-a}\phi ^q \\
+c\int_M\left| Rm\right| ^p\left( f+1\right) ^{-a}\phi
^q+ca^{2p}\int_M\left( f+1\right) ^{-a+p},
\end{gather*}
which by (\ref{h}) yields

\begin{equation}
\left( a-c\right) \int_M\left| Rm\right| ^p\left( f+1\right) ^{-a}\phi
^q\leq c\int_M\left| \nabla Rm\right| ^2\left| Rm\right| ^{p-3}\left(
f+1\right) ^{-a}\phi ^q+C,  \label{k}
\end{equation}
where 
\begin{eqnarray*}
C &:&=C_1+C_2+C_3, \\
C_3 &:&=ca^{2p}\int_M\left( f+1\right) ^{-a+p}.
\end{eqnarray*}
Finally, let us use (\ref{f2}) to check that 
\begin{gather}
2\int_M\left| \nabla Rm\right| ^2\left| Rm\right| ^{p-3}\left( f+1\right)
^{-a}\phi ^q\leq \int_M\left( \Delta \left| Rm\right| ^2\right) \left|
Rm\right| ^{p-3}\left( f+1\right) ^{-a}\phi ^q  \label{l} \\
-\int_M\left( \nabla f\cdot \nabla \left| Rm\right| ^2\right) \left|
Rm\right| ^{p-3}\left( f+1\right) ^{-a}\phi ^q+c\int_M\left| Rm\right|
^p\left( f+1\right) ^{-a}\phi ^q.  \notag
\end{gather}
The first term in the formula above is 
\begin{gather*}
\int_M\left( \Delta \left| Rm\right| ^2\right) \left| Rm\right| ^{p-3}\left(
f+1\right) ^{-a}\phi ^q \\
=-\int_M\left( \nabla \left| Rm\right| ^2\cdot \nabla \left| Rm\right|
^{p-3}\right) \left( f+1\right) ^{-a}\phi ^q \\
+a\int_M\left( \nabla \left| Rm\right| ^2\cdot \nabla f\right) \left|
Rm\right| ^{p-3}\left( f+1\right) ^{-a-1}\phi ^q \\
-q\int_M\left( \nabla \left| Rm\right| ^2\cdot \nabla \phi \right) \left|
Rm\right| ^{p-3}\left( f+1\right) ^{-a}\phi ^{q-1}.
\end{gather*}
Choosing $p\geq 3$ guarantees that

\begin{equation*}
-\int_M\left( \nabla \left| Rm\right| ^2\cdot \nabla \left| Rm\right|
^{p-3}\right) \left( f+1\right) ^{-a}\phi ^q\leq 0.
\end{equation*}
On the other hand, 
\begin{gather*}
a\int_M\left( \nabla \left| Rm\right| ^2\cdot \nabla f\right) \left|
Rm\right| ^{p-3}\left( f+1\right) ^{-a-1}\phi ^q \\
\leq ac\int_M\left| \nabla Rm\right| \left| Rm\right| ^{p-2}\left(
f+1\right) ^{-a}\phi ^q \\
\leq \frac 14\int_M\left| \nabla Rm\right| ^2\left| Rm\right| ^{p-3}\left(
f+1\right) ^{-a}\phi ^q+ca^2\int_M\left| Rm\right| ^{p-1}\left( f+1\right)
^{-a}\phi ^q.
\end{gather*}
Similarly we find:

\begin{gather*}
-q\int_M\left( \nabla \left| Rm\right| ^2\cdot \nabla \phi \right) \left|
Rm\right| ^{p-3}\left( f+1\right) ^{-a}\phi ^{q-1} \\
\leq \frac 14\int_M\left| \nabla Rm\right| ^2\left| Rm\right| ^{p-3}\left(
f+1\right) ^{-a}\phi ^q \\
+c\int_M\left| Rm\right| ^{p-1}\left( f+1\right) ^{-a}\phi ^{q-2}.
\end{gather*}
Notice moreover that 
\begin{gather*}
-\int_M\left( \nabla f\cdot \nabla \left| Rm\right| ^2\right) \left|
Rm\right| ^{p-3}\left( f+1\right) ^{-a}\phi ^q \\
\leq 2\int_M\left| \nabla Rm\right| \left| \nabla f\right| \left| Rm\right|
^{p-2}\left( f+1\right) ^{-a}\phi ^q \\
\leq \frac 12\int_M\left| \nabla Rm\right| ^2\left| Rm\right| ^{p-3}\left(
f+1\right) ^{-a}\phi ^q+c\int_M\left| Rm\right| ^{p-1}\left| \nabla f\right|
^2\left( f+1\right) ^{-a}\phi ^q.
\end{gather*}
Using these estimates in (\ref{l}) we get: 
\begin{gather*}
\int_M\left| \nabla Rm\right| ^2\left| Rm\right| ^{p-3}\left( f+1\right)
^{-a}\phi ^q\leq c\int_M\left| Rm\right| ^p\left( f+1\right) ^{-a}\phi ^q \\
+ca^2\int_M\left| Rm\right| ^{p-1}\left( f+1\right) ^{-a+1}\phi ^{q-2}.
\end{gather*}
From interpolation as in formula (\ref{g}), for $\varepsilon =\frac 1{a^2},$
it follows: 
\begin{eqnarray*}
&&ca^2\int_M\left| Rm\right| ^{p-1}\left( f+1\right) ^{-a+1}\phi ^{q-2} \\
&\leq &c\int_M\left| Rm\right| ^p\left( f+1\right) ^{-a}\phi
^q+ca^{2p}\int_M\left( f+1\right) ^{-a+p}\phi ^{q-2p}.
\end{eqnarray*}
Therefore, we have proved: 
\begin{equation*}
\int_M\left| \nabla Rm\right| ^2\left| Rm\right| ^{p-3}\left( f+1\right)
^{-a}\phi ^q\leq c\int_M\left| Rm\right| ^p\left( f+1\right) ^{-a}\phi
^q+ca^{2p}\int_M\left( f+1\right) ^{-a+p}.
\end{equation*}
This, by (\ref{k}) implies that 
\begin{equation*}
\left( a-c\right) \int_M\left| Rm\right| ^p\left( f+1\right) ^{-a}\phi
^q\leq C.
\end{equation*}
Recall that $c$ is a constant depending only on $n,p$ and $K$ while $a$ is a
sufficiently large arbitrary number. This shows that there exists $a$
depending on $n,p$ and $K$ such that 
\begin{equation*}
\int_M\left| Rm\right| ^p\left( f+1\right) ^{-a}\phi ^q\leq C.
\end{equation*}
To conclude the proof of the Lemma notice that if 
\begin{equation*}
\rho \left( x\right) \leq \frac 12r,
\end{equation*}
then 
\begin{equation*}
\phi \left( x\right) \geq \frac 3{16},
\end{equation*}
hence this shows that 
\begin{equation*}
\int_{D\left( \frac 12r\right) }\left| Rm\right| ^p\left( f+1\right)
^{-a}\leq C.
\end{equation*}
Since $C$ is independent of $r,$ making $r\rightarrow \infty $ we get 
\begin{equation*}
\int_M\left| Rm\right| ^p\left( f+1\right) ^{-a}\leq C.
\end{equation*}
Moreover, from $C=C_1+C_2+C_3$ and the expressions for these constants we
see that in fact we have the estimate 
\begin{equation*}
\int_M\left| Rm\right| ^p\left( f+1\right) ^{-a}\leq c\int_{D\left(
r_0\right) }\left| Rm\right| ^p\left( f+1\right) ^{-a}+c\int_M\left(
f+1\right) ^{-a+p},
\end{equation*}
where $r_0$ is a fixed number, depending on $n,p$ and $K.$

This proves the Lemma in the case $p\geq 3$. We can interpolate as in (\ref
{g}) to get the claim for any $p\geq 2.$ The claim that $\int_{B_x\left(
1\right) }\left| Rm\right| ^p\leq C\left( r\left( x\right) +1\right) ^{2a}$
follows immediately from the asymptotics of $f,$ see above.  $\blacksquare $

We are now ready to finish the proof of the Theorem. We will be brief here,
since this part of the proof is standard.

From (\ref{f2}) we infer that 
\begin{eqnarray*}
\Delta \left| Rm\right| ^2 &\geq &2\left| \nabla Rm\right| ^2+\nabla f\cdot
\nabla \left| Rm\right| ^2-c\left| Rm\right| ^3 \\
&\geq &-c\left( \left| Rm\right| +\left| \nabla f\right| ^2\right) \left|
Rm\right| ^2.
\end{eqnarray*}
Therefore, if we denote 
\begin{equation*}
u:=c\left( \left| Rm\right| +\left| \nabla f\right| ^2\right) ,
\end{equation*}
then 
\begin{equation*}
\Delta \left| Rm\right| ^2\geq -u\left| Rm\right| ^2.
\end{equation*}
Since we assumed the Ricci curvature is bounded below, there exists a
uniform bound on the Sobolev constant of the ball $B_x\left( 1\right) .$
More exactly, for any $\varphi $ with support in $B_x\left( 1\right) $ we
have: 
\begin{equation*}
\left( \int_{B_x\left( 1\right) }\varphi ^{\frac{2n}{n-2}}\right) ^{\frac{n-2%
}n}\leq C_S\int_{B_x\left( 1\right) }\left( \left| \nabla \varphi \right|
^2+\varphi ^2\right) .
\end{equation*}
Indeed, we can use the Sobolev inequality in \cite{SC} and notice that we
can control the volume of $B_x\left( 1\right) $ from below uniformly in $x$
since $M$ satisfies a log Sobolev inequality and the Ricci curvature is
bounded, see \cite{CN}. The constant $C_S$ depends only on $n$, the Ricci
curvature lower bound on $M$ and on $\inf_{x\in M}vol\left( B_x\left(
1\right) \right) .$ As proved by Carillo and Ni in \cite{CN}, $vol\left(
B_x\left( 1\right) \right) \geq k>0,$ where $k$ depends on $K=\sup_M\left|
Rc\right| $ and on Perelman's $\mu $ invariant. We remark that with our
normalization of $f,$ the $\mu $ invariant can be computed as $\mu =\log
\left( \left( 4\pi \right) ^{\frac n2}\int_Me^{-f}\right) <\infty .$

Then the standard Moser iteration, see \cite{L}, implies that 
\begin{equation*}
\left| Rm\right| ^2\left( x\right) \leq A\int_{B_x\left( 1\right) }\left|
Rm\right| ^2,
\end{equation*}
where 
\begin{equation*}
A:=\bar{C}\left( \int_{B_x\left( 1\right) }u^n+1\right) ,
\end{equation*}
for a constant $\bar{C}$ depending only on $n$ and $C_S.$ Since we showed in
Lemma 2 that $\int_{B_x\left( 1\right) }u^p$ grows at most polynomially in $%
r\left( x\right) ,$ the theorem follows from here.$\blacksquare $

\section{Gap and compactness theorems}

In this section we prove the gap theorem and the compactness theorem of
shrinking Ricci solitons, based on the estimates proved above. Since
everywhere in this section the Ricci curvature is bounded, we can apply
Theorem 1 to see that $\left| Rm\right| \left( x\right) \leq C\left(
1+r\left( x\right) \right) ^{2a}.$ We follow a similar argument as in
Theorem 1, this time using the weight $e^{-f}$ and paying more attention to
the dependence on the Ricci curvature bound. The computations will be
simpler, and we will use many times the identities 
\begin{eqnarray*}
\nabla _l\left( R_{ijkl}e^{-f}\right) &=&0, \\
\nabla _j\left( R_{ij}e^{-f}\right) &=&0.
\end{eqnarray*}
Clearly, we do not need to use a cut-off here, since all the curvature terms
will be integrable with respect to $e^{-f},$ by Theorem 1. As in the
previous section, we take 
\begin{equation*}
K=\sup_M\left| Rc\right| .
\end{equation*}
We also assume that $K>0,$ since otherwise there is nothing to prove. We use
the notation $\Delta _f=\Delta -\nabla f\cdot \nabla $ and note that $\Delta
_f$ is self adjoint with respect to the weighted volume $e^{-f}dv.$

Since 
\begin{equation*}
\Delta _f\left( f\right) =\Delta f-\left| \nabla f\right| ^2=\frac n2-f,
\end{equation*}
it follows that 
\begin{gather}
\int_M\left| Rm\right| ^p\left( f-\frac n2\right) e^{-f}=-\int_M\left|
Rm\right| ^p\Delta _f\left( f\right) e^{-f}=\int_M\nabla f\cdot \nabla
\left| Rm\right| ^pe^{-f}  \label{1} \\
=p\int_M\nabla _hR_{ijkl}R_{ijkl}f_h\left| Rm\right|
^{p-2}e^{-f}=2p\int_M\nabla _lR_{ijkh}R_{ijkl}f_h\left| Rm\right|
^{p-2}e^{-f}  \notag \\
=-2p\int_MR_{ijkh}f_{hl}R_{ijkl}\left| Rm\right|
^{p-2}e^{-f}-2p\int_MR_{ijkh}f_hR_{ijkl}\nabla _l\left( \left| Rm\right|
^{p-2}\right) e^{-f}.  \notag
\end{gather}
As in Theorem 1, we take $p\geq 3.$ Using the soliton equation and the Ricci
curvature bound $\left| Rc\right| \leq K,$ we find that 
\begin{eqnarray*}
R_{ijkh}f_{hl}R_{ijkl} &=&\frac 12\left| Rm\right| ^2-R_{ijkh}R_{ijkl}R_{hl}
\\
&\geq &\left( \frac 12-K\right) \left| Rm\right| ^2.
\end{eqnarray*}
We can also estimate 
\begin{gather*}
-2p\int_MR_{ijkh}f_hR_{ijkl}\nabla _l\left( \left| Rm\right| ^{p-2}\right)
e^{-f}=-4p\int_M\left( \nabla _jR_{ik}\right) R_{ijkl}\nabla _l\left( \left|
Rm\right| ^{p-2}\right) e^{-f} \\
\leq 4p\left( p-2\right) \int_M\left| \nabla Rc\right| \left| \nabla
Rm\right| \left| Rm\right| ^{p-2}e^{-f}.
\end{gather*}
Then (\ref{1}) shows that 
\begin{equation}
\int_M\left( f-\frac n2+p\left( 1-2K\right) \right) \left| Rm\right|
^pe^{-f}\leq 4p^2\int_M\left| \nabla Rc\right| \left| \nabla Rm\right|
\left| Rm\right| ^{p-2}e^{-f}.  \label{2}
\end{equation}
We estimate the right hand side of (\ref{2}) as follows: 
\begin{gather}
2\int_M\left| \nabla Rc\right| \left| \nabla Rm\right| \left| Rm\right|
^{p-2}e^{-f}\leq \frac 1{pK}\int_M\left| \nabla Rc\right| ^2\left| Rm\right|
^{p-1}e^{-f}  \label{2b} \\
+pK\int_M\left| \nabla Rm\right| ^2\left| Rm\right| ^{p-3}e^{-f}.  \notag
\end{gather}
Furthermore, we have 
\begin{eqnarray*}
\Delta _f\left| Rc\right| ^2 &=&2\left| \nabla Rc\right| ^2+2\left|
Rc\right| ^2-4R_{ikjh}R_{ij}R_{kh} \\
&\geq &2\left| \nabla Rc\right| ^2-4K^2\left| Rm\right| .
\end{eqnarray*}
Consequently, 
\begin{gather*}
\int_M\left| \nabla Rc\right| ^2\left| Rm\right| ^{p-1}e^{-f}\leq \frac
12\int_M\left( \Delta _f\left| Rc\right| ^2\right) \left| Rm\right|
^{p-1}e^{-f}+2K^2\int_M\left| Rm\right| ^pe^{-f} \\
=-\frac 12\int_M\nabla \left| Rc\right| ^2\cdot \nabla \left| Rm\right|
^{p-1}e^{-f}+2K^2\int_M\left| Rm\right| ^pe^{-f} \\
\leq \left( p-1\right) K\int_M\left| \nabla Rc\right| \left| \nabla
Rm\right| \left| Rm\right| ^{p-2}e^{-f}+2K^2\int_M\left| Rm\right| ^pe^{-f}.
\end{gather*}
Using this in (\ref{2b}) we get 
\begin{equation*}
\int_M\left| \nabla Rc\right| \left| \nabla Rm\right| \left| Rm\right|
^{p-2}e^{-f}\leq \frac{2K}p\int_M\left| Rm\right| ^pe^{-f}+pK\int_M\left|
\nabla Rm\right| ^2\left| Rm\right| ^{p-3}e^{-f},
\end{equation*}
which, after plugging into (\ref{2}), yields 
\begin{equation}
\int_M\left( f-\frac n2+p\left( 1-10K\right) \right) \left| Rm\right|
^pe^{-f}\leq 4p^3K\int_M\left| \nabla Rm\right| ^2\left| Rm\right|
^{p-3}e^{-f}.  \label{3}
\end{equation}
Finally, we have: 
\begin{equation*}
2\int_M\left| \nabla Rm\right| ^2\left| Rm\right| ^{p-3}e^{-f}\leq
\int_M\left( \Delta _f\left| Rm\right| ^2\right) \left| Rm\right|
^{p-3}e^{-f}+20\int_M\left| Rm\right| ^pe^{-f}.
\end{equation*}
Indeed, to see this one only has to check the details of the proof of (\ref
{f2}). Integrating by parts and using that $p\geq 3$ we find 
\begin{equation*}
\int_M\left( \Delta _f\left| Rm\right| ^2\right) \left| Rm\right|
^{p-3}e^{-f}=-\int_M\left( \nabla \left| Rm\right| ^2\cdot \nabla \left|
Rm\right| ^{p-3}\right) e^{-f}\leq 0.
\end{equation*}
Therefore, from (\ref{3}) we conclude that for $p\geq 3$ we have: 
\begin{equation}
\int_M\left( f-\frac n2+p\left( 1-50p^2K\right) \right) \left| Rm\right|
^pe^{-f}\leq 0.  \label{4}
\end{equation}
We are now ready to prove the Corollaries.

\textit{Proof of Corollary 1.}

Let us take $p=n.$ We check from (\ref{4}) that if $K\leq \frac 1{100n^2}$
then 
\begin{equation*}
\int_Mf\left| Rm\right| ^ne^{-f}\leq 0.
\end{equation*}
Recall that $f\,$is normalized such that $\left| \nabla f\right| ^2+R=f,$
and in particular, since any gradient shrinker has $R\geq 0$ (see \cite{C,
Ca}) it follows that $f\geq 0.$ Thus the above inequality implies that $M$
is flat i.e. $\left( M,g,f\right) $ is the Gaussian soliton $\left( %
\mathbb{R}^n,dx^2,\frac 14\left| x\right| ^2\right) .\blacksquare $

\textit{Proof of Corollary 2.}

We take $p=\frac n2$ in (\ref{4}) to see that 
\begin{equation}
\int_M\left( f-7n^3K\right) \left| Rm\right| ^{\frac n2}e^{-f}\leq 0.
\label{5}
\end{equation}
We fix $x_0$ a point where $f$ achieves its minimum on $M$. Then we have
(see \cite{CZ, HM}) 
\begin{equation}
\frac 14\left[ \left( d\left( x_0,x\right) -5n\right) _{+}\right] ^2\leq
f\left( x\right) \leq \frac 14\left( d\left( x_0,x\right) +\sqrt{2n}\right)
^2,  \label{7}
\end{equation}
where $a_{+}:=\max \left\{ 0,a\right\} .$

Let us set $r_0:=6n+\sqrt{28n^3K}.$ Using this in (\ref{5}) we see that: 
\begin{gather*}
\int_{M\backslash B_{x_0}\left( r_0\right) }\left| Rm\right| ^{\frac
n2}e^{-f}\leq \int_{f\geq 7n^3K+1}\left| Rm\right| ^{\frac n2}e^{-f} \\
\leq \int_{f\geq 7n^3K+1}\left( f-7n^3K\right) \left| Rm\right| ^{\frac
n2}e^{-f}\leq 7n^3K\int_{f\leq 7n^3K+1}\left| Rm\right| ^{\frac n2}e^{-f} \\
\leq 7n^3K\int_{B_{x_0}\left( r_0\right) }\left| Rm\right| ^{\frac n2}\leq
7n^3KL.
\end{gather*}
Using again (\ref{7}) shows that for any $r>0$ 
\begin{eqnarray*}
\int_{B_{x_0}\left( r\right) }\left| Rm\right| ^{\frac n2} &\leq &E\left(
r\right) ,\;\;\text{for} \\
E\left( r\right) &:&=7n^3KLe^{\frac 14\left( r+\sqrt{2n}\right) ^2}.
\end{eqnarray*}

This proves Corollary 2. $\blacksquare $

\bigskip

{\tiny DEPARTMENT OF MATHEMATICS, COLUMBIA UNIVERSITY, NEW YORK, NY 10027}%
\newline
{\small E-mail address: omuntean@math.columbia.edu}

\bigskip

{\tiny DEPARTMENT OF MATHEMATICS, COLUMBIA UNIVERSITY, NEW YORK, NY 10027}%
\newline
{\small E-mail address: mtwang@math.columbia.edu}

\end{document}